\documentclass{jnmp01b}

\usepackage{amsmath}

\numberwithin{equation}{section}

 \newtheorem{theorem}{Theorem}[section]
 \newtheorem{prop}{Proposition}[section]
\theoremstyle{definition}
 \newtheorem{note}{Note}[section]
 \newtheorem{example}{Example}[section]
 \newtheorem*{comment}{Comment}

\newcommand{\Cn}{\mathop{\rm Cn}\nolimits}
\newcommand{\Sn}{\mathop{\rm Sn}\nolimits}

\makeatletter
\DeclareRobustCommand{\primfrac}[1]{%
  \PackageWarning{amsmath}{%
Foreign command \@backslashchar#1; %
\protect\frac\space or \protect\genfrac\space should be used instead%
  }
  \global\@xp\let\csname#1\@xp\endcsname\csname @@#1\endcsname
  \csname#1\endcsname
}
\makeatother

\begin{document}

\renewcommand{\evenhead}{Yu.A.\ Il'in}
\renewcommand{\oddhead}{Asymptotic Integration of Nonlinear Systems}


\thispagestyle{empty}

\begin{flushleft}
\footnotesize \sf
Journal of Nonlinear Mathematical Physics \qquad 2000, V.7, N~2,
\pageref{firstpage}--\pageref{lastpage}.
\hfill {\sc Article}
\end{flushleft}

\vspace{-5mm}

\copyrightnote{2000}{Yu.A.\ Il'in}

\Name{Asymptotic Integration of Nonlinear Systems
of Differential Equations whose Phase Portrait
       is  Foliated on Invariant Tori}

\label{firstpage}

\Author{Yuri A.\ IL'IN}

\Adress{Faculty of Mathematics and Mechanics, St.-Petersburg State University,\\
Bibliotechnaja 2, St. Petergof, St.-Petersburg, 198804, Russia \\
E-mail: iljin@paloma.spbu.ru}

\Date{Received February 8, 2000; Accepted March 12, 2000}

\begin{abstract}
\noindent
We consider the class of autonomous systems $\dot x=f(x)$, where
 $x \in {\bf R}^{2n}$, $f~\in~C^1({\bf R}^{2n})$ whose
 phase portrait is a Cartesian product of $n$  two-dimensional {\em centres}.
We also consider  perturbations of this system, namely
$\dot x=f(x)+g(t,x)$, where
$g \in C^1({\bf R}\times{\bf R}^{2n})$ and $g$ is
asymptotically small, that is $g\Rightarrow 0$
as $t\to +\infty$ uniformly with respect to $x$. The rate of
decrease of $g$ is assumed to be $t^{-p}$ where $p>1$.
We prove under this conditions the existence of
bounded solutions of the perturbed system and discuss their convergence to
solutions of the unperturbed system. This convergence depends on  $p$.
Moreover, we show that the original unperturbed system may be reduced
to the form  $\dot r=0$, $\dot\theta=A(r)$,
and taking $r\in {\bf R}^m_{+}$,
 $\theta\in {\bf T}^n$, where ${\bf T}^n$ denotes the $n$-dimensional
 torus, we investigate the more general case of systems whose phase portrait
is foliated on invariant tori. We notice that integrable
 Hamiltonian systems are of the same nature. We give also
 several examples, showing that the conditions of our theorems cannot be
 improved.
\end{abstract}


\section {Introduction}

\paragraph*{A.}
In \cite{Il} the author  investigated
 the asymptotic behaviour
 of solutions of an ${\bf R}^2 $--system
\[
\dot x=f(x)+g(t,x),
\]
where the unperturbed system $\dot x=f(x)$ has a phase portrait like a
{\it linear centre} and the perturbation $g(t,x)\Rightarrow 0$
as $t \to +\infty$.
In this paper we give natural generalizations of this problem.

\paragraph*{B.}
First of all we  consider an ${\bf R}^{2n}$--system
\begin{equation}
\left\{
   \begin{array}{rcl}
       \dot x_1 & = & f_1 (x_1)\\
               &\dots&  \\
       \dot x_n & = & f_n (x_n),\\
   \end{array}
\right.           
\end{equation}
where $x_k \in {\bf R}^2$ and $f_k \in {\rm C}^1({\bf R}^2 , {\bf R}^1)$
for all $k=1,\dots,n$. We also assume that the phase portrait of each
subsystem $\dot x_k =f_k (x_k)$ is similar to a {\it linear centre},
i.e. origin is a unique equilibrium point and the other trajectories are closed
curves surrounding the origin. The phase portrait of the whole system (1.1)
resembles  the Cartesian product
\begin{equation}
({\bf R}^1_+ \times {\bf S}^1) \times \dots \times
                   ({\bf R}^1_+ \times {\bf S}^1)=
        {\bf R}^n_+ \times {\bf T}^n,  
\end{equation}
where ${\bf R}^1_+=[0, +\infty)$ and ${\bf T}^n$ denotes the $n$--dimensional
torus.

 The first example of such a system is a linear system $\dot x=Ax$
 where all eigenvalues of matrix $A$ are pure imaginary and have only
 prime Jordan blocks. This system can be reduced by a linear
 transformation  to a system of $n$ independent mathematical pendulums
 \[
\left\{
   \begin{array}{ccccc}
       y''_1&+&a_1 y_1& = & 0\\
               &\dots&&&  \\
       y''_n&+&a_n y_n& = & 0\\
   \end{array}
\right.
\]
where all $a_k > 0$.
The next more general example is a system of $n$ independent Duffing equations
\[
\left\{
   \begin{array}{ccccc}
       y''_1&+&a_1 y_1& = & h_1(y_1,y'_1)\\
               &\dots&&&  \\
       y''_n&+&a_n y_n& = & h_n(y_n,y'_n)\\
   \end{array}
\right.
\]
where every function $h_k$ does not destroy the centre defined by the
linear part in each equation. We often meet the same systems in  different
applications of mathematics, mechanics, physics and engineering
sciences.
There are other  important examples of similar systems,
but we do not mention all of them here. The only
exception is a Hamiltonian integrable system mentioned below.

We rewrite  system (1.1) in brief form
\begin{equation}
          \dot x=f(x), \tag{1.1}
\end{equation}
where $x=(x_1,\dots,x_n)^{\top}$, $f=(f_1,\dots,f_n)^{\top}$.
And we consider a perturbation of (1.1)
\begin{equation}
             \dot x=f(x)+g(t,x), 
\end{equation}
where $g\in {\rm C}^{0,1}_{t,x}({\bf R}^1_t \times
{\bf R}^{2n}_x , {\bf R}^1)$
and function $g$ satisfies the main assumption
\begin{itemize}
  \item[({\bf A1})]
    {\em $ t^p \|g(t,x)\|\Rightarrow 0 $
       uniformly with respect to $x$ on  compact subsets
       of  ${\bf R}^{2n}$ as $t \to \infty$. }
\end{itemize}
An equivalent formulation of (A1) is:
\begin{itemize}
  \item[({\bf A2})]
    {\em there are  positive continuous scalar functions $\alpha(s)$
    and $\beta(t)$ such that $\alpha$ is increasing, $\beta$ is decreasing,
    $\lim_{t \to \infty}\beta(t)=0$ and  }
\[
     t^p \|g(t,x)\| \le  \alpha(\|x\|)\beta(t).
\]
\end{itemize}
The proof of this equivalence is standard enough and we omit it here.
The following theorems are proved for system (1.3).

{\bf Theorem 4.2.}   If $p>1$ then system (1.3) has  bounded
                    solutions for $t \ge 0$.

{\bf Theorem 4.4.}   If $p>1$ then any bounded solution of (1.3)
    approaches (orbitally) a certain $n$-dimensional invariant torus
    of the unperturbed system (1.1).

{\bf Theorem 5.2.}   If $p>2$ then any bounded solution of (1.3)
    approaches  a certain solution of unperturbed system (1.1).

 It will be shown that for $p=1$ and $p=2$  these theorems will
 no longer be valid.

\paragraph*{C.}
In Section~3 we  construct for the system (1.1)
action-angle variables (``quasi-polar'' coordinates) in which
(1.1) has the form  $ \dot r = 0$,
    $ \dot \theta  =  A(r) $
where $r=(r_1,\dots,r_n)^{\top}$ are ``quasi-polar'' radii and
$\theta=(\theta_1,\dots,\theta_n)^{\top}$ are ``quasi-polar'' angles.

 It gives us possibility  to consider instead of (1.1) the more general
 system  in the form
 \begin{equation}
 \left\{
     \begin{array}{rcl}
     \dot r & = & 0\\
     \dot \theta & = & A(r),\\
     \end{array}
\right.            
\end{equation}
where $r \in {\bf R}^m_+$, $\theta \in {\bf T}^n$,
$ A \in {\rm C}^1({\bf R}^m_+ , {\bf R}^n_+ \setminus \{0\})$
(i.e. each component $A_k(r)>0$). The phase portrait of (1.4)
resembles the Cartesian product of a more complicated form
 than (1.2), namely
\[
({\bf R}^{m_1}_+ \times {\bf T}^n_1) \times \dots \times
                   ({\bf R}^{m_d}_+ \times {\bf T}^{n_d})=
        {\bf R}^m_+ \times {\bf T}^n
\]
 where $m_1+\dots+m_d=m$ and $n_1+\dots+n_d=n.$
The perturbations of (1.4) need  to be taken in the form
 \begin{equation}
 \left\{
     \begin{array}{rcl}
     \dot r & = & P(t,r,\theta)\\
     \dot \theta & = & A(r)+Q(t,r,\theta)\\
     \end{array}
\right.            
\end{equation}
where $P \in {\rm C}^{0,1,1}_{t,r,\theta}
({\bf R}^1 \times {\bf R}^m_+ \times {\bf T}^n , {\bf R}^m)$,
$Q \in {\rm C}^{0,1,1}_{t,r,\theta}
({\bf R}^1 \times {\bf R}^m_+ \times {\bf T}^n , {\bf R}^n)$
(Let us recall that $\theta \in {\bf T}^n$ means indeed that
$\theta \in {\bf R}^n$ and functions $P$ and $Q$ are periodic in
each component of the vector $\theta $.)

We shall suppose that $P$ and $Q$ satisfy the following main assumption
\begin{itemize}
  \item[({\bf A3})]
     {\em $t^p\|P(t,r,\theta)\|, t^p\|Q(t,r,\theta)\| \Rightarrow 0$
  uniformly with respect to $r$ on compact subsets of  ${\bf R}^m_+$
  and for all $\theta \in {\bf T}^n$  as $t \to \infty$.}
\end{itemize}
An equivalent formulation of (A3) is:
\begin{itemize}
  \item[({\bf A4})]
   {\em there are  positive continuous scalar functions $\alpha(s)$
   and $\beta(t)$ such that $\alpha$ is increasing,
   $\beta$ is decreasing,
   $\lim_{t \to \infty}\beta(t)=0$ and }
  \[
  t^p \|P(t,r,\theta)\|,t^p \|Q(t,r,\theta)\| \le  \alpha(\|r\|)\beta(t).
  \]
\end{itemize}
We prove  the following theorems for (1.5).

{\bf Theorem 4.1.}   If $p>1$ then system (1.5) has solutions
            bounded with respect to $r$-coordinates for $t \ge 0$.

{\bf Theorem 4.3.}   If $p>1$ then any bounded solution of (1.5)
    approaches (orbitally) a certain  invariant torus
    of the unperturbed system (1.4).

{\bf Theorem 5.1.}   If $p>2$ then any bounded solution of (1.5)
    converges to a certain solution of the unperturbed system (1.4).

 We note at the end of this section that
 problems of asymptotic integration are
 classical in the theory of differential equations. If we eliminate
 the $\theta$--equation and
 the variable $\theta$ from (1.5) we obtain a well--investigated in
 the literature problem of the perturbations of zero system
 (see~\cite{W1}--~\cite{W3},
 monographs~\cite{Ha} and \cite{Y1} give the general theory).
 From this point of view Theorems
 4.1--4.4 are quite natural and  expected statements.
 On the contrary,
 Theorems 5.1--5.2 represent more delicate  results.
We consider the construction of the ``quasi-polar'' variables and proof of
Theorems 5.1--5.2 as one of the most important parts of our paper.

\section{The existence of a positively-definite integral}

The aim of this section is to prove the following statement.
\begin{theorem}
  Let $\dot x_k=f_k(x_k)$ be a two-dimensional subsystem of (1.1).
  Then it has a differentiable positively--definite integral $U_k(x_k)$
  on ${\bf R}^2$.
\end{theorem}
 This is a known
 result. But because we shall need further
  some detailed properties of $U_k$ we give an explicit proof
 of this theorem.

\begin{proof}
For convenience we omit the index $k$ and will write
simply $x$ and $f$ instead of $x_k$ and $f_k$.
We belief that no confusion
with the system (1.1) can arise. Let $x=(z_1, z_2)^\top$  and
$f=(h_1, h_2)^\top$. The system $\dot x=f(x)$ can be rewritten in the form
\begin{equation}
  \dot z_1=h_1(z_1,z_2),\quad \dot z_2=h_2(z_1,z_2).    
\end{equation}
Consider the system ``orthogonal'' to  (2.1)
\begin{equation}
  \dot z_1 = - h_2 (z_1, z_2),\quad \dot z_2 = h_1 (z_1, z_2 ).
\end{equation}
 It is evident that the right-hand vector fields of (2.1) and (2.2)
 are orthogonal at each point. This implies that (2.2) has also
 a unique equilibrium point at the origin and other trajectories of (2.2)
 cross orthogonally all cycles of (2.1) at a unique point. For if any
 trajectory of (2.2) intersects some cycle of (2.1) at more than
 one point then there is a point on this cycle where right-hand vectors
 of (2.1) and (2.2) will tangent each other.

Let $\xi(t)$ be a certain nonzero solution of (2.2) and $I=(a,b)$ denote
its maximal interval of existence. Let $\Gamma$ denote the trajectory
of $\xi(t)$. Without loss of generality we can assume that $\Gamma$ crosses
the cycles of (2.1) from inside to outside with
 increase of $t$ (otherwise
 one needs to consider
 the system  $\dot z_1=h_2,\; \dot z_2=-h_1$ instead of (2.2)).
 It follows from Poincar\'e--Bendixon
 theory that $\xi(t) \to 0$ as $t \to a_+$
 and $\|\xi(t)\| \to \infty$ as $t \to b_-$. Indeed,  system (2.2) cannot
 have  closed curves (cycles) because then they must intersect the
 cycles of (2.1) at more than one point.
  Hence the $\alpha$-limit set of $\xi(t)$ is not
 empty and may consist of only origin.
 By the same reasons the $\omega$-limit set
 of $\xi(t)$ must be empty, i.e. $\|\xi(t)\| \to \infty$ as $t \to b_-$.
   It is clear that  when  time $t$ increases from $a$
    to $b$ the solution $\xi(t)$ crosses
 each cycle of (2.1). We have also that $a=-\infty$ because $\xi(t)$
 remains bounded as $t \to a_+$.

 Let $\varphi(t,x^0)$ denote the solution of (2.1) with initial condition
 $\varphi(0,x^0)=x^0$. We define the mapping $\chi : {\bf R}^2 \setminus \{0\}
 \to \Gamma$ by the formula
 \[
   \chi(x)=\{\varphi(t,x): t \in {\bf R}\} \cap \Gamma =
              \varphi({\bf R},x) \cap \Gamma.
 \]
 It is clear that $\chi \in {\rm C}^1$ because the cycle $\varphi({\bf R},x)$
 intersects  $\Gamma$ transversally. Let $\xi^{-1}:\Gamma \to I$ be
 the mapping which to every $x\in \Gamma$ associates the time $t$ such that
 $x=\xi(t)$. Since $\xi(t)$ is differentiable and $ \xi'(t) \ne 0$
(it is not equilibrium point) we have $\xi^{-1}\in {\rm
C}^1(\Gamma, I)$.
Finally we define on ${\bf R}^2\setminus\{0\}$ the function $U_1(x)$ by the
equality
\[
   U_1(x)=(\xi^{-1}\circ \chi)(x).
\]
We point out   some properties of $U_1$.
\begin{itemize}
  \item[{\bf 1})]
     $U_1 \in {\rm C}^1({\bf R}^2 \setminus \{0\} , I)$,
  \item[{\bf 2})]
     $U_1(\xi(t))=t$,
  \item[{\bf 3})]
     $U_1(\varphi(t,x)) \equiv U_1(x)$,  thus $U_1$ is an
     integral for  system (2.1) on   ${\bf R}^2 \setminus \{0\}$,
  \item[{\bf 4})]
     equation $U_1(x)=r$ defines a unique cycle of (2.1)\\
     (indeed, curve $\xi(t)$
     intersects all cycles of (2.1) with increasing $t$ at only
     one point; hence  $U_1$ monotonically increases along $\Gamma$;
     therefore each level curve of $U_{1}(x)$ consists of a unique cycle).
  \item[{\bf 5})]
     $\lim_{x \to 0} U_1(x)=-\infty$. (Indeed, $x \to 0 \Rightarrow
     \chi(x) \to 0 \Rightarrow (\xi^{-1}\circ \chi)(x) \to a=-\infty $).
\end{itemize}

 By setting $U_2(x)=\exp(U_1(x))$ we
 construct a continuous positively-definite
integral for (2.1) on whole ${\bf R}^2$  (we assume naturally $U_2(0)=0$).
It is clear that $U_2({\bf R}^2)=[0,e^b).$
 Since differentiability of $U_2$ can be broken only at $x=0$ let us
consider the behaviour of the
 gradient $DU_2(x)$ as $x \to 0$. There are two possibilities.
\begin{itemize}
  \item[{\bf 1})]
    $\limsup\limits_{x \to 0} \|DU_2(x)\| < \infty$.
  \item[{\bf 2})]
    $\limsup\limits_{x \to 0} \|DU_2(x)\| = \infty$.
\end {itemize}
In the first case we define $U(x)=U_2^2(x)$. We thus have
\(
\lim_{x \to 0}DU(x)\\=2\lim_{x \to 0}U_2(x)DU_2(x) = 0
\)
and $U(x)$ is therefore the required integral for (2.1).

In the second case we will construct the integral $U(x)$ in the form
$(\mu \circ U_2)(x)$ where a scalar function $\mu$ is chosen to smooth out
the discontinuity of $DU_2$ at origin. To do this we first define the
function $\nu(s)=\sup\limits_{x \in (U_2)^{-1}(s)} \|DU_2(x)\|$ for
$s \in (0,e^b)$.
We have  $\nu \in {\rm C}(0,e^b)$,
since the sets $(U_2)^{-1}(s)$ are compact and
 depend continuously on $s$.
The following properties are evident:
\begin{itemize}
  \item[{\bf 1})]
     $\nu(s)\ge 0$.
  \item[{\bf 2})]
     $\lim_{s \to +0}\nu(s)= \infty$.
  \item[{\bf 3})]
     $\nu(U_2(x))\ge \|DU_2(x)\|$.
\end {itemize}
Finally we define $\mu(s)=\int_0^s (1+\nu(u))^{-2}du$. It is easy to see
that $\mu \in {\rm C}^1[0,e^b)$ and $\mu(0)=0$. We set
$U(x)=\mu(U_2(x))$ and verify that $U(x)$ is the required integral for (2.1).
All we need is to check that $DU(x)$ exists at $x=0$.
But it is fulfilled, because
\begin{gather*}
\lim_{x \to 0} \|DU(x)\|=\lim_{x \to 0} \|\mu'(U_2(x))DU_2(x)\|
\\
\qquad=\lim_{x \to 0} \left\|\frac{DU_2(x)}{(1+\nu(U_2(x)))^2}\right\|
\stackrel{(property \,3)}{\le}
\lim_{x \to 0} \left(\frac{1}{1+\|DU_2(x)\|}\right)=0.
\end{gather*}
We point out some properties of  $U$.
\begin{itemize}
  \item[{\bf 1})]
     $U \in {\rm C}^1({\bf R}^2 , [0,b^*))$, where $b^*=e^{2b}$
         or $\mu(e^b)$.
  \item[{\bf 2})]
     $U(x)=0 \Leftrightarrow x=0$, i.e. $U$ is positively definite.
  \item[{\bf 3})]
     $U$ monotonically increases along  $\Gamma : \quad
      \partial U/\partial\Gamma>0$.
  \item[{\bf 4})]
     each level curve $U(x)=r$ defines a unique cycle of (2.1).
\end{itemize}
The theorem is proved.
\end{proof}

\begin{note}
For simplicity we will further assume that $b^*=+\infty$
because it is easy to construct a ${\rm C}^1$--function
$\gamma : [0,b^*)\to[0,+\infty)$ and consider the new integral
$\gamma(U(x))$ instead of $U(x)$.
\end{note}

\section{The action-angle  or quasi-polar coordinates}

In this section we construct  new coordinates for more convenient
notation of  systems (1.1) and (1.3).
This procedure is a standard one for Hamiltonian systems, see for
example the textbook~\cite{Ta}.
But we would notice that our coordinates are not identical to the hamiltonian
action-angle coordinates. We introduce them in a more complicated and
sly way, so that the change of variables does not break the  asymptotic
properties of the perturbations and solutions.

To each ${\bf R}^2$--component
$x_k$ of vector $x$ we put in correspondence two  scalar coordinates
$r_k$ and $\theta_k$.
We define
the first ``quasi-polar'' coordinate (the ``radius'') $r_k$
by $r_k=U_k(x_k)$.
In accordance to Section~2, the
equation $U_k(x_k)=r_k$ for each $r_{k}>0$
defines a unique cycle of $\dot x_k=f_k(x_k)$.
Moreover the function $r_k=U_k(x_k)$ is invertible on $\Gamma$,
since $U_k$ monotonically increases along  $\Gamma$.
Let $v_k(r_k)$ denote this inversion. We have
\begin{itemize}
  \item[{\bf 1})]
     $v_k \in {\rm C}^1({\bf R}_+ , \Gamma)$,
  \item[{\bf 2})]
     $v_k(r_k)=\Gamma \cap \{x_k : U_k(x_k)=r_k \}$,
  \item[{\bf 3})]
     $U_k(v_k(r_k))=r_k$.
\end{itemize}
We denote the period of any cycle  $U_k(x_k)=r_k$ by $T_k(r_k)$. It is
equal to the time for the next intersection with $\Gamma$
of the solution passing at $t=0$ through the point $x_k=v_k(r_k)$ .
Since $\Gamma$ is a tranversal, we have $T_k \in {\rm C}^1$.
We now prove the following
\begin{prop}
  The value $1/T_k(r_k)$ is well-defined on  ${\bf R}_+$.
\end{prop}
\begin{proof}
Obviously, $T_k(r_k)\ne 0$ for $r_k\ne 0$. It  remains
to prove that $T_k(0)\ne 0$. Let us  use the notations from the proof
of Theorem 2.1. Changing in (2.1) coordinates $z_1$ and $z_2$ into the
standard polar coordinates $z_1=\rho \cos \alpha$,
$z_2=\rho \sin \alpha$, we get for $\alpha$
\[
  \dot\alpha=\rho^{-1}(\cos\alpha h_2(\rho \cos \alpha, \rho \sin \alpha)-
   \sin\alpha h_1(\rho \cos \alpha, \rho \sin \alpha))
    =\rho^{-1}F(\rho, \alpha).
\]
Obviously, $F \in {\rm C}^1({\bf R}_+\times {\bf R})$
and $F(0,\alpha)\equiv 0$.
Then Hadamard's lemma implies $F(\rho, \alpha)\equiv \rho G(\rho, \alpha)$
where $G \in {\rm C}({\bf R}_+\times {\bf R})$. Hence,
$ \dot\alpha=G(\rho, \alpha)$. Let $(\rho(t), \alpha(t))$ denote the solution
passing through the initial point $(\rho_0, 0)$ at $t=0$. Let $T(\rho_0)$
be the period of it. It is clear that $\alpha(T(\rho_0))=2\pi$. On other
hand, $\alpha(T(\rho_0))=\int_0^{T(\rho_0)} G(\rho(t), \alpha(t))\,dt$.
For any $\rho_0 \in [0,1]$ the solution $(\rho(t), \alpha(t))$
belongs to a certain
compact set, therefore $G(\rho(t), \alpha(t))$ is uniformly bounded with
respect to $\rho_0$. If $T(\rho_0)\to 0$ as $\rho_0\to 0$, then
\[
  2\pi=\lim_{T(\rho_0)\to 0} \alpha (T(\rho_0))
      =\lim_{T(\rho_0)\to 0}\int_0^{T(\rho_0)} G(\rho(t), \alpha(t))\,dt
      =0.
\]
This contradiction proves the proposition.
\end{proof}

Change  time in the system $\dot x_k=f_k(x_k)$ by formula
\[
  t=T_k(r_k)\tau.
\]
Since $\dot r_k=0$ along the solutions we have $dt=T_k(r_k)d\tau$ and
\begin{equation}
   \frac{dx_k}{dt}=T_k(r_k)f_k(x_k).  
\end{equation}
All solutions of (3.1) have the common period $1$, thus the centre defined
by (3.1) is an isochronous one. Let $\psi_k(\tau,x_k)$ be the solution of (3.1)
satisfying the initial condition $\psi_k(0,x_k)=x_k$.
We define the second ``quasi--polar''
coordinate (the ``angle'') $\theta_k$ in such a way that
$\psi_k(-\theta_k,x_k)\in \Gamma$. I.e. $-\theta_k$ is the time of the first
intersection of solution $\psi_k(\tau,x_k)$ with $\Gamma$.
It is clear
that $\theta_k$ depends  smoothly on $x_k$. Thus the general formula
for the transformation to the new coordinates is
\begin{equation}
  x_k=\psi_k(\theta_k,v_k(r_k)). 
\end{equation}

Let us investigate some properties of the ``quasi--polar'' variables. For any
function $h(x,y)$ the notations $D_1 h$ and $D_2 h$ denote the partial
derivatives of $h$ with respect to the first and second argument
correspondingly.
Let $J_k(r_k,\theta_k)$ be the Jacobian matrix of the change of variables (3.2).
We have
\[
  J_k(r_k,\theta_k)=\frac{D\psi_k(\theta_k,v_k(r_k))}{D(r_k,\theta_k)} =
  \big[D_2\psi_k(\theta_k,v_k(r_k)) \, v'_k(r_k),
   \quad D_1\psi_k(\theta_k,v_k(r_k))\big].
\]
Since $\psi_k(\theta_k+1,v_k(r_k))\equiv \psi_k(\theta_k,v_k(r_k))$,
the matrix
$D_2\psi_k(\theta_k,v_k(r_k))$ is  1--periodic  with respect
to $\theta_k$. On the other hand $D_2\psi_k$ is the
derivative of the  solution
$\psi_k$ with respect to initial values, i.e. it is a fundamental matrix
of a so-called linear system in variation. Hence $D_2\psi_k$ is nonsingular
and $D_2\psi_k(0,v_k(r_k))=E$. It follows from periodicity that
the values 
\[
\|D_2\psi_k(\theta_k,v_k(r_k))\|,\ \ 
\|(D_2\psi_k(\theta_k,v_k(r_k)))^{-1}\|,\ \
\det D_2\psi_k(\theta_k,v_k(r_k))
\]
and $\det(D_2\psi_k(\theta_k,v_k(r_k)))^{-1}$ are bounded away from zero
and infinity for each fixed $r_k$ and $\forall \theta_k \in {\bf R}$.

Since $\psi_k(\theta_k,v_k(r_k))$ is a solution of (3.1), we have
\[
D_1\psi_k(\theta_k,v_k(r_k))=T_k(r_k)f_k(\psi_k(\theta_k,v_k(r_k))).
\]
This yields $J_k(r_k,0)=\big[E\cdot v'_k(r_k),\, T_k(r_k)f_k(v_k(r_k))\big]$.
But derivative
$v'_k$ is a tangent vector to  $\Gamma$ at the point $v_k(r_k)$,
because $v_k(r_k)$ is some parametrization of curve $\Gamma$.
Hence $v_k'(r_k)$ is to be orthogonal to the vector $f_k(v_k(r_k))$.
Thus matrix
$J_k(r_k,0)$ is non singular $\forall r_k > 0$.

It is easy to check the following statement.
Let $\dot y=F(y)$ be an arbitrary
$C^1$--smooth system and $\varphi(t,y)$ be the solution of it with
initial condition $\varphi(0,y)=y$, then
$F(\varphi(t,y))=D_2\varphi(t,y)F(y)$.

This implies
$T_k(r_k)f_k(\psi_k(\theta_k,v_k(r_k)))=D_2\psi_k(\theta_k,v_k(r_k))
                  T_k(r_k)f_k(v_k(r_k)) $
and therefore
\begin{align*}
  J_k(r_k,\theta_k)&=\big[(D_2\psi_k)\, v'_k,\quad (D_2\psi_k)T_k\,f_k\big]=
  (D_2\psi_k)  \big[ v'_k,\quad T_k\,f_k\big]
\\
  &=D_2\psi_k(\theta_k,v_k(r_k)) J_k(r_k,0).
\end{align*}

Each matrix in the last expression is nonsingular and possesses
``good'' properties. So this equality guarantees
that the Jacobian matrix $J_k$ is nonsingular
for fixed $r_k >0$ and
arbitrary $\theta_k \in {\bf R}$.
Moreover, $\det J_k$, $\det J_k^{-1}$, $\|J_k\|$ and $\|J_k^{-1}\|$ are
bounded away from zero and infinity. Similar properties of boundedness
will also be fulfilled if we  take $r_k \in K$ where $K$ is some
compact set of  ${\bf R}_+$ and $\theta_k$ as before
will be an arbitrary one.

In the new variables  system (3.1) will have the  following form
\[
  \frac{dr_k}{d\tau}=0, \quad \frac{d\theta_k}{d\tau}=1.
\]
The original system $\dot x_k = f_k(x_k)$ is transformed into
\[
  \dot r_k =0, \quad \dot \theta_k =1/T_k(r_k)
\]
(replacing $dt = T_k(r_k)d\tau$).

Finally the full system (1.1) is transformed into
\begin{equation}
  \dot r=0, \quad  \dot \theta = A(r), 
\end{equation}
where $r=(r_1,\dots,r_n)^{\top}$, $\theta=(\theta_1,\dots,\theta_n)^{\top}$,
$A(r)=(A_1(r_1),\dots,A_n(r_n))^{\top}$,  $A_k(r_k)=1/T_k(r_k)$.

Changing the variables in the perturbed system (1.3), we obtain
\begin{equation}
  {\dot r \choose \dot\theta}=
     {0\choose A(r)}+[J(r,\theta)]^{-1}g(t,\psi(\theta, v(r)))=
  {0\choose A(r)}+{P(t,r,\theta)\choose Q(t,r,\theta)} 
\end{equation}
(we use natural notations $J={\rm diag}(J_1,\dots,J_n)$,
$\psi=(\psi_1,\dots, \psi_n)^{\top}$,\\ $v(r)=(v_1(r_1),\dots,v_n(r_n))^{\top}$).
The functions $P$ and $Q$ possess the following properties:
\begin{itemize}
  \item[{\bf 1})]
      they are in $C^{0,1,1}_{t,r,\theta}$,
  \item[{\bf 2})]
      they are 1--periodic with respect to each components of $\theta$,
  \item[{\bf 3})]
      if $g$ satisfies the assumptions (A1) (or (A2))
      then $P$ and $Q$ satisfy (A3) (or (A4)).
\end{itemize}
Indeed, let $K \subset {\bf R}^n_+$ be a compact set and $r\in K$. There exist
$M>0$ such that  $\|(J(r,\theta))^{-1}\| \le M$ for all $r \in K$ and
$\theta \in {\bf R}^n$. Hence
\[
 t^p\|{P(t,r,\theta) \choose Q(t,r,\theta)}\|
   \le Mt^p\|g(t,\psi(\theta,v(r)))\| \Rightarrow 0
     \quad {\rm as} \quad t \to +\infty.
\]
(we note that if $r\in K$ then cycle $\psi(\theta,v(r))$ does not
 leave some
 compact set in ${\bf R}^{2n}$ for all $\theta \in {\bf R}^n$).

 Thus, we have shown that systems (1.1) and (1.3) may be reduced to systems
(1.4) and (1.5) while preserving assumption (A3) (or (A4)).

\begin{note}
 We want to underline that it is not so difficult
 to reduce formally system
 (1.3) to (1.5)  (see theory of Hamiltonian system in~\cite{Ta}).
 The problem is to keep the asymptotic properties (A1-A2) as well as
the boundedness
 and covergence of solutions under the inverse transformation.
 For this we use the trick with isochronicity. Without it the Jacobian matrix
 may be unbounded and this manner of reasoning fails.
\end{note}

\section{The existence of bounded solutions}

Let us consider the system (1.5).
\begin{theorem}
  Let $P$ satisfy assumption (A3) or (A4) with $p>1$.
  Then (1.5) has solutions bounded with respect to $r$-coordinates.
\end{theorem}
\begin{proof}
We shall find a comparison equation for $\|r\|$. We have
\[
  \frac{d}{dt}\|r\|^2=\frac{d}{dt}(r,r)=2(r,\dot r)=2(r,P(t,r,\theta)),
\]
here $(\cdot,\cdot)$ denotes inner product. This gives for $\|r\|\ne 0$
\begin{equation}
  \frac{d}{dt}\|r\|\le \|P(t,r,\theta)\|\stackrel{(A4)}{\le}
   \alpha(\|r\|)\beta(t)/t^p \le \alpha(\|r\|)/t^p, 
\end{equation}
(we can assume that $\beta(t)\le 1$ for sufficiently large $t$, because
$\beta(t)\to 0$). Let us consider the comparison equation
\begin{equation}
  \frac{d}{dt}\rho = \alpha(\rho)/t^p.  
\end{equation}
It is well--known (see~\cite{Ha}) that every function $r(t)$ satisfying
inequality (4.1) does not exceed the solution $\rho(t)$ of (4.2) with initial
condition $\rho(t_0)=\|r(t_0)\|$. Thus the existence of bounded solutions
of (4.2)
will imply the same for (1.5).

Let $G(\rho)=\int_1^\rho 1/\alpha(s)\,ds$ .
Since $G'(\rho)>0$, function  $G(\rho)$ is monotone increasing
and invertible. The solution of (4.2) with initial
condition $\rho(t_0)=\rho_0$  has the form
\begin{equation}
  G(\rho(t))-G(\rho_0)=(p-1)(t_0^{1-p}-t^{1-p}). 
\end{equation}

There are two possibilities.

{\bf A}. $\lim_{\rho\to\infty}G(\rho)=+\infty$.
Then all solutions of (4.2) must be bounded. For if some solution
$\rho(t)\to +\infty$ as $t \to +\infty$, then the left--hand side of (4.3)
tends to $\infty$ whereas the right--hand side tends to
$(p-1)t_0^{1-p}<\infty$. We got a contradiction.

{\bf B}. $\lim_{\rho\to\infty}G(\rho)=G^*<\infty$. Fix $\rho_0 >0$ and take
$t_0$ so large that $G^*-G(\rho_0)>(p-1)t_0^{1-p}$.
Then the solution of (4.2) satisfying the initial  condition
$\rho(t_0)=\rho_0$ is bounded. Otherwise, tending $t\to+\infty$
in (4.3) and assuming that $\rho(t)\to +\infty$ we get the impossible
equality $G^*-G(\rho_0)=(p-1)t_0^{1-p}$.

Thus the equation (4.2) has  bounded solutions in  both cases.
So the theorem is proved.
\end{proof}

\begin{note}
The case {\bf A} occurs when $P$ satisfies
the well--known inequality (see~\cite{Ha})
\[
 \|P(t,r,\theta)\|\le h(\|r\|)\gamma(t),
\]
where $\int^\infty \gamma(t)\,dt<\infty$  and $\int^\infty 1/h(s)\,ds=\infty$.
In this case all solutions of (1.5) are bounded with respect
to $r$. Conversely, there are unbounded solutions in the case {\bf B}
 as the next example shows.
\end{note}

\begin{example}
Consider the system $\dot r=0$, $\dot\theta=1$ and its
perturbation $\dot r=r^2/t^2$, $\dot\theta=1$. The perturbed system has
unbounded extendable to infinity solution $r=t$. All solutions lying inside
the sector $0<r<t$ are bounded and those lying outside it are unbounded
(moreover, they are unextendable to infinity).
\end{example}

\begin{note}
The condition $p>1$ in Theorem 4.1 cannot be generalized
to include the case $p=1$. See Example 4.2 below.
\end{note}

Now we give reformulations of Theorem 4.1 and Note 2.1 for system
(1.3).
\begin{theorem}
   Suppose that perturbation $g$ in (1.3) satisfies (A1) or (A2) with $p>1$ .
   Then (1.3) has  bounded solutions.
\end{theorem}
\begin{proof}
As it follows from Section~3 the system (1.3) written in
``quasi-polar'' coordinates satisfies all conditions of Theorem 4.1.
\end{proof}

\begin{note}
All solutions of (1.3) are bounded
if $\|g(t,x)\|\le h(\|x\|)\gamma(t)$, where
$\int^\infty \gamma(t)\,dt<\infty$  and $\int^\infty 1/h(s)\,ds=\infty$.
\end{note}

\begin{example}
The present example shows that condition $p>1$
cannot be generalized to include the case
 $p=1$. We take the standard  two-dimensional centre
\begin{equation}
 \dot x_1=-x_2, \quad \dot x_2=x_1,  
\end{equation}
and perturb it in the following way
\[
\dot x_1=-x_2+\frac{\cos t}{t\ln t}, \quad \dot x_2=x_1+\frac{\sin t}{t\ln t},
 \quad (t>1).
\]
It is evident that $\cos t/(t\ln t)$ and $\sin t/(t\ln t)$ satisfy (A2)
with $p=1$, $\alpha(s)=1$, $\beta(t)=1/\ln t$. But all solutions of the
perturbed system have the form
\[
x_1=(C_1+\ln (\ln t))\cos t - C_2\sin t, \quad
x_2=(C_1+\ln (\ln t))\cos t + C_2\sin t,
\]
and are unbounded.
\end{example}

Now we return to system (1.5) and show that under the assumptions of
Theorem 4.1 every bounded solution approaches a certain invariant torus
of (1.4).
\begin{theorem}
 Let $(r(t),\theta(t))$ be a bounded solution of (1.5)
 with respect to $r$-coordinates and
\begin{equation}
   t^p\|P(t,r(t),\theta(t))\| \to 0 \quad {\rm as}
     \quad t \to 0 \tag{$\tilde A3$}
\end{equation}
with $p>1$. There exists $r^*$ such that
\[
   \lim_{t\to+\infty}r(t)=r^*.
\]
\end{theorem}

\begin{comment}
It means that solution $(r(t),\theta(t))$ approaches
the invariant torus $r=r^*$ of (1.4).
\end{comment}

\begin{note}
Evidently, if (A3) or (A4) holds then every bounded
solution of (1.5) satisfies $(\tilde A3)$.
\end{note}

\begin{proof}
We have
\[
 r(t)=r(t_0) + \int_{t_0}^t \dot r(s)\,ds=
      r(t_0) + \int_{t_0}^t P(s,r(s),\theta(s))\,ds.
\]
Assumption $(\tilde A3)$ implies
that $\int_{t_0}^\infty \|P(s,r(s),\theta(s))\|\,ds$
converges. Hence the limit
\[
 \lim_{t\to+\infty} r(t)=r(t_0) + \int_{t_0}^\infty P(s,r(s),\theta(s))\,ds=r^*
\]
exists. The theorem is proved.
\end{proof}

\begin{note}
The condition $p>1$ cannot be generalized
to include the case $p=1$, see Example~4.4.
\end{note}

\begin{note}
The next example demonstrates that assumptions of Theorem~4.3
don't imply even orbital convergence of bounded solutions of (1.5) to certain
solution of (1.4). However this convergence
occurs if $\theta \in {\bf R}^1$
 (see~\cite{Il}), see also Note 5.2.
\end{note}

\begin{example}
Consider the following system
\begin{equation}
 \dot r_1=0 , \quad \dot r_2=0 , \quad  \dot\theta_1=r_1 ,
  \quad \dot\theta_2=r_2 , 
\end{equation}
and its perturbation
\begin{equation}
\dot{\bar r}_1=1/t^2 , \quad  \dot{\bar r}_2=2/t^2 ,
 \quad \dot{\bar\theta}_1={\bar r}_1 ,
 \quad \dot{\bar\theta}_2={\bar r}_2 , 
\end{equation}
Solutions of (4.5) are
\begin{equation}
 r_1=r^0_1 ,  \quad r_2=r^0_2 ,  \quad \theta_1=r_1^0 t+\theta_1^0 ,
  \quad      \theta_2=r_2^0 t+\theta_2^0 , 
\end{equation}
where $r^0_1$, $r^0_2$, $\theta_1^0$, $\theta_2^0$ are arbitrary constants
(initial values). Solutions of (4.6) are
\begin{equation}
{\bar r}_1={\bar r}^0_1-1/t ,  \quad {\bar r}_2={\bar r}^0_2-2/t ,  \quad
{\bar\theta}_1={\bar r}_1^0 t-\ln t +{\bar\theta}_1^0 ,  \quad
{\bar\theta}_2={\bar r}_2^0 t-2\ln t +{\bar\theta}_2^0 . 
\end{equation}
All solutions of (4.6) are bounded with respect to $\bar r_1$ and $\bar r_2$.
It is clear that $\bar r_1 \to \bar r_1^0$
and $\bar r_2 \to \bar r_2^0$ as
$t \to \infty$. The orbital convergence of (4.8) to (4.7) means that
$\bar r_k^0= r_k^0 ,\,\, k=1,2$ and there is an infinitely monotonically
increasing function $s(t)$ such that
\[
 \lim_{t\to\infty}(\bar \theta_1(t)-\theta_1(s(t)))=
 \lim_{t\to\infty}(\bar \theta_2(t)-\theta_2(s(t)))=0,
\]
or
\begin{gather}
 \lim_{t\to\infty}(r_1^0 t-\ln t +\bar\theta_1^0 -r_1^0s(t) - \theta_1^0)=0,
\\
 \lim_{t\to\infty}(r_2^0 t-2\ln t +\bar\theta_2^0 -r_2^0s(t) - \theta_2^0)=0.
\end{gather}
From (4.9) we have
$s(t)=t-(\ln t -\bar\theta_1^0 + \theta_1^0)/r_1^0 + o(1)$.
Substituting this in (4.10) we obtain
\[
 \lim_{t\to\infty}((\frac{r_2^0}{r_1^0}-2)\ln t + \bar \theta_2^0
   -\theta_2^0 +(\bar \theta_1^0-\theta_1^0)/r_1^0 + o(1))=0.
\]
However, this is possible if and only  if $r_2^0=2r_1^0$,
$\bar \theta_1^0= \theta_1^0$, $\bar \theta_2^0= \theta_2^0$.
Consequently there are solutions
of (4.6) which don't orbitally converge to any solution of (4.5).
\end{example}

Now we reformulate Theorem 4.3 for  system (1.3).
\begin{theorem}
  Let $x(t)$ be a bounded solution of (1.3) and
\[
    t^p\|g(t,x(t))\| \to 0 \quad {\rm as} \quad t \to +\infty
\]
 with $p>1$. Then $x(t)$ approaches a certain invariant torus of
 unperturbed system (1.1).
\end{theorem}
\begin{proof}
It follows immediately from Theorem 4.3 and Section~3.
\end{proof}

\begin{note}
For Theorem 4.4 Notes 4.4--4.6 are  also valid.
\end{note}

\begin{example}
The next example demonstrates that assumption $p>1$
cannot be generalized to include the case $p=1$.
Let's consider again the two-dimensional linear centre (4.4)
with the following perturbation
\begin{equation}
\dot x_1=-x_2+\frac{\cos(\ln(\ln t))}{t\ln t} x_{1}, \quad
\dot x_2=x_1+\frac{\sin(\ln(\ln t))}{t\ln t} x_{2},
 \quad (t>1). 
\end{equation}
Obviously every bounded solution
of (4.11) satisfies the conditions of Theorem 4.4 with $p=1$. But (4.11)
has the bounded solutions
\[
 x_1=r_0\cos t \exp(\sin(\ln(\ln t))), \quad
 x_2=r_0\sin t \exp(\sin(\ln(\ln t))),
\]
which oscillate between the circles with radii $r_{0}e$ and
$r_{0}e^{-1}$. Hence no one of them can approach any circle
$ x_1=r\cos t , \;  x_2=r\sin t $, which is just the
solution of the unperturbed system.
\end{example}

\section{The convergence of bounded solutions}

We begin with the main statement for system (1.5).
\begin{theorem}
 Let (A3) or (A4) hold with $p>2$. Then each bounded solution of (1.5)
 converges to a certain solution of the  unperturbed system (1.4).
\end{theorem}
\begin{proof}
Let $(r(t),\theta(t))$ denote any bounded solution of (1.5).
Theorem 4.1 implies the convergence $r(t)\to r^*$ as $t\to +\infty$.
Every solution of (1.4)
lying on the invariant torus $r=r^*$ has the form $(r^*, A(r^*)t+\theta_0)$,
where $\theta_0$ is the initial value. It is necessary and sufficient to prove
the existence of a $\theta_0$ such that
\[
  \lim_{t\to+\infty}(\theta(t)-A(r^*)t-\theta_0)=0.
\]
Thus it is sufficient to prove that the limit value
$\theta_0=\lim_{t\to+\infty}(\theta(t)-A(r^*)t)$ exists.
Since
\begin{gather*}
 \theta(t)-A(r^*)t=\theta(1)-A(r^*)+
        \int_1^t (\dot\theta(s)-A(r^*))\,ds
\\
\qquad  =\theta(1)-A(r^*)+\int_1^t (A(r(s))+Q(s,r(s),\theta(s))-A(r^*))\,ds ,
\end{gather*}
it is sufficient to prove that the obtained integral converges.
By (A4) we have
\begin{equation}
 \|{P\choose Q}(s,r(s),\theta(s))\| \le \alpha(r(s))\beta(s)/s^p. 
\end{equation}
There is a constant $M$ such that $\alpha(r(s))\beta(s)\le M$ for $s>1$.
Hence
\[
 \|\int_1^\infty Q(s,r(s),\theta(s))\,ds\|\le
   M\int_1^\infty 1/s^p\,ds=M/(p-1).
\]
Next we need to estimate $\|A(r(s))-A(r^*)\|$. First of all we notice that
\begin{gather}
 \|r(s)-r^*\|=\|\int_{+\infty}^s \dot r(u)\,du\| \le
  \int_s^{+\infty} \|P(u,r(u),\theta(u))\,du\|
\notag\\
  \qquad \le\int_s^{+\infty} \alpha(r(u))\beta(u)/u^p\,du \le
  M\int_s^{+\infty}1/u^p\,du = M s^{1-p}/(p-1).  
\end{gather}
Since $A\in {\rm C}^1$, it satisfies the Lipschitz inequality with some
constant $L$ (on such a compact set $K$ from ${\bf R}^m_+$
that $r(s)\in K$ for $s\ge 1$ )
\begin{equation}
  \|A(r)-A(r^*)\| \le L\|r-r^*\|. 
\end{equation}
 It  follows from (5.2) and (5.3) that $\|A(r(s))-A(r^*)\|=O(s^{1-p})$.
 Inequality $p>2$ implies convergence of
  $\int_s^{+\infty}\|A(r(s))-A(r^*)\|\,ds$ . Hence limit
\[
  \lim_{t\to+\infty}(\theta(t)-A(r^*)t)=
  \theta(1)-A(r^*)+\int_1^{\infty} (A(r(s))+Q(s,r(s),\theta(s))-A(r^*))\,ds
\]
 exists. This proves the theorem.
\end{proof}

 Now we give the formulation of the analogous theorem for (1.3).
 \begin{theorem}
 Let (A1) or (A2) hold with $p>2$. Then each bounded solution of (1.3)
 converges to a certain solution of the unperturbed system (1.1).
\end{theorem}

\begin{note}
The condition $p>2$ cannot be generalized
to include the case $p=2$ as the following
example shows.
\end{note}

\begin{example}
Consider the system
\begin{equation}
 \dot x_1=-x_2^3, \quad \dot x_2=x_1^3, 
\end{equation}
 and its perturbation
\begin{equation}
 \dot x_1=-x_2^3+\frac{1+\ln t}{(t\ln t)^2}x_1,
 \quad \dot x_2=x_1^3+\frac{1+\ln t}{(t\ln t)^2}x_2. 
\end{equation}
The cycles of (5.4) are $x_1^4+x_2^4=C,$ $C\in {\bf R}_+$.
Let $(\Cn t, \Sn t)$ denote the solution of (5.4) passing through
point (1,0) as $t=0$. Changing the variables by the formulas
 $x_1=r\Cn\theta,\, x_2=r\Sn\theta$  (action-angle coordinates)
we transform (5.4) into
\begin{equation}
        \dot r=0, \quad
        \dot \theta=r^2, 
\end{equation}
and (5.5) into
\begin{equation}
        \dot r=r\frac{1+\ln t}{(t\ln t)^2}, \quad
        \dot\theta =r^2. 
\end{equation}
The solutions of (5.6) are $r(t)=r_0,$  $\theta(t)=\theta_0+r_0^2t$
and those of (5.7) are
$\bar r(t)=r_0\exp[-1/(t\ln t)],$
$\bar\theta(t)=\theta_0+r_0^2\int_2^t \exp[-2/(s\ln s)]\,ds.$
Obviously, $\lim \bar r(t)=r_0.$ But
$ \lim\limits_{t\to+\infty}(\theta(t)-A(r_0)t)$ from proof of Theorem 5.1
does not exist. Indeed,
\begin{gather*}
  \lim_{t\to+\infty}(\theta(t)-A(r_0)t)=
  \lim_{t\to+\infty}(\theta(t)-r_0^2t)=
  \lim_{t\to+\infty}r_0^2 \big(\int_2^t \exp[-2/(s\ln s)]\,ds-t\big)
\\
  \qquad =\lim_{t\to+\infty}r_0^2\big(\int_2^t (\exp[-2/(s\ln
s)]-1)\,ds-2\big)=-\infty,
\end{gather*}
because obtained integral diverges. This follows from the asymptotic
equivalence
\[
   \exp[-2/(s\ln s)]-1 \sim -2/(s\ln s)
\]
and divergence of $\int_2^{\infty} 2/(s\ln s)\,ds.$
\end{example}

\begin{note}[about convergency in Section~4] 
We  wouls like to make an interesting
remark about the possibility to decrease the lower bound in the inequality
$p>2$ in Theorems 5.1 and 5.2.
Suppose $A(r)\equiv A_0$ is a constant. For system (1.1) it is equivalent
to $T_k(r_k)\equiv T_k.$ It means that each subsystem
$\dot x_k=f_k(x_k)$ in (1.1), $k=1, \dots ,n,$  defines an isochronous
centre. For (1.5) we can also talk about a property of isochronicity.
Then the statements
of Theorems 5.1--5.2 are valid for $p>1.$ Indeed, we have  in this case
\[
  \lim_{t\to+\infty}(\theta(t)-A(r^*)t)=
  \theta(1)-A(r^*)+\int_1^{\infty} Q(s,r(s),\theta(s))\,ds ,
\]
where convergence of this integral follows from (5.1).
\end{note}

The property of isochronicity in Hamiltonian systems
(and not only for such systems) is discussed in numerous papers.
We refer only to~\cite{CMV}, where further references can be found.
Only for asymptotic perturbations of
such systems, our Note 5.2 is valid.

We note that the main reason of appearance of the condition $p>2$
is just the necessity to estimate  $A(r(t))-A(r^*)$.
Hence we may suppose
that $Q$ satisfies (A3) with $p>1$ (but $P$ with $p>2$ as before)
and Theorem 5.1 will also be valid.

Finally we point out a possibility of the situation for system (1.1)
when some subsystems of it
are isochronous centres and all other  not. It is not difficult to see
that for ``isochronous'' angles $\theta_k$ the orbital convergence
({\em Note:} even usual convergence)  is realized in Theorem 4.4.

\subsection*{Acknowledgements}

I am very grateful to the Royal Swedish Academy of Sciences for their
financial
support of this work. I would also like to thank G. S\"oderbacka
for carefully reading this paper and his valuable remarks.

\label{lastpage}

\end{document}